\documentclass[11pt]{article}
\usepackage{amsfonts}
\usepackage{color}
\usepackage{graphics}
\usepackage{indentfirst}

\usepackage{cite}
\usepackage{latexsym}
\usepackage{amsmath}
\usepackage{amssymb}
\usepackage[dvips]{epsfig}
\usepackage{amscd}

 \usepackage{color}
 \usepackage{indentfirst}
\usepackage{ntheorem}
\theoremstyle{plain}
\theoremseparator{.}
\newtheorem{theorem}{Theorem}[section]
\newtheorem{remark}{Remark}[section]

\newtheorem{lemma}[theorem]{Lemma}

\newcommand\thmref[1]{Theorem~\ref{#1}}
\newcommand\lemref[1]{Lemma~\ref{#1}}

\newcommand{\ti}{\tilde}

\def\pf{{\it Proof.}  }

\newcommand{\thatsall}{\hfill$\Box$}
\newcommand{\bi}{\bibitem}

\newcommand{\bt}{\begin{theorem}}
\newcommand{\bl}{\begin{lemma}}
\newcommand{\el}{\end{lemma}}
\newcommand{\et}{\end{theorem}}

\renewcommand{\b}{\beta  }

\newcommand{\te}{\theta}

\newcommand{\al}{\alpha}

\newcommand{\ve}{\varepsilon}
\newcommand{\la}{\label}

\newcommand{\ka}{\kappa}

\newcommand{\bn}{\begin{eqnarray}}
\newcommand{\en}{\end{eqnarray}}
\newcommand{\bnn}{\begin{eqnarray*}}
\newcommand{\enn}{\end{eqnarray*}}

\newcommand{\bnnn}{\begin{eqnarray*}}
\newcommand{\ennn}{\end{eqnarray*}}
\newcommand{\ben}{\begin{enumerate}}
\newcommand{\een}{\end{enumerate}}
\newcommand{\ba}{\begin{aligned}}
\newcommand{\ea}{\end{aligned}}
\newcommand{\be}{\begin{equation}}
\newcommand{\ee}{\end{equation}}

\def\norm[#1]#2{\|#2\|_{#1}}

\def\xix{\int_0^1}

\makeatletter
\@addtoreset{equation}{section}
\makeatother

\makeatletter
\@addtoreset{equation}{section}
\makeatother

\title{ Global Strong Solutions to Compressible Navier-Stokes System  with
Degenerate Heat Conductivity  and Density-Depending Viscosity \thanks{ Partially supported by NNSFC   11671027 and 11471321.}
}

\author{Bin Huang,  Xiaoding Shi, Ying Sun \thanks{   Email addresses: abinhuang@gmail.com, abinhuang36@163.com (B. Huang), shixd@mail.buct.edu.cn (X. Shi), 1913349041@qq.com (Y. Sun)}  \\[3mm]     Department of Mathematics, Faculty of Science, \\Beijing University of Chemical Technology, \\ Beijing  100029, P. R. China }
\date{ }

\begin{document}
\maketitle
\begin{abstract}
We consider the compressible Navier-Stokes system where  the viscosity depends on density   and the heat conductivity is  proportional to a positive power of  the temperature
under stress-free and thermally insulated boundary conditions. Under the same   conditions   on  the  initial data  as those of the constant viscosity and   heat conductivity case ([Kazhikhov-Shelukhin. J. Appl. Math. Mech. 41
(1977)], we obtain the existence and uniqueness of global strong solutions.
Our result can be regarded as a natural generalization of the Kazhikhov's theory  for the constant heat conductivity case to the degenerate and nonlinear one under stress-free and thermally insulated boundary conditions.  \end{abstract}

\noindent{\bf Keywords:}  Compressible Navier-Stokes system; Density-depending viscosity;
Degenerate heat conductivity;  Stress-free

\section{Introduction}

The compressible Navier-Stokes system which  describes  the  one-dimensional motion of a
viscous heat-conducting   perfect polytropic gas is written
  in the Lagrange variables in the following form (see \cite{4,23}):
  \be\la{1.1}
v_t=u_{x},
\ee
\be\la{1.2}
u_{t}+P_{x}=\left(\mu\frac{u_{x}}{v}\right)_{x},\ee
 \be\la{1.3}\left(e+\frac{1}{2}u^2\right)_t+\left(Pu\right)_x=\left(\frac{\kappa\te_x+\mu uu_x}{v}\right)_{x}, \ee
where $t>0$ is time, $x\in \Omega= (0,1)$ denotes the
Lagrange mass coordinate,  and the unknown functions $v>0, u$ and $P$ are,  respectively, the specific volume of the gas, fluid velocity,  and  pressure. In this paper, we concentrate on ideal polytropic gas, that is, $P$ and $e$ satisfy \be \la{1.4}   P =R \theta/{v},\quad e=c_v\theta +\mbox{const},
\ee
where  both specific gas constant  $R$ and   heat
capacity at constant volume $c_v $ are   positive constants. We also assume that $\mu$ and $\ka$ satisfy \be\la{1.5}   \mu=\ti\mu \left(1+v^{-\al}\right), \quad \ka=\ti\ka \te^\beta, \ee
with constants $\ti\mu,\ti\ka>0$ and $\al,\beta\ge 0.$
The system \eqref{1.1}-\eqref{1.5} is supplemented
with the initial  conditions
\be\la{1.6}(v,u,\te)(x,0)=(v_0,u_0,\te_0)(x),  \quad  x\in (0,1), \ee
under stress-free and thermally insulated boundary  conditions\be\la{1.7} \left(\frac{\mu}{v}u_x-P\right)(0,t)=\left(\frac{\mu}{v}u_x-P\right)(1,t)=0,\quad\te_x(0,t)=\te_x(1,t)=0,\ee
and the initial data \eqref{1.6} should be compatible with the boundary conditions \eqref{1.7}.

 The boundary conditions \eqref{1.7} describe the expansion of a finite mass of gas into vacuum.
One also considers other kind of boundary conditions
\be\ba \la{1.77} u(0,t) =u(1,t)=0,\quad\te_x(0,t)=\te_x(1,t)=0,\ea\ee
which mean that the gas is confined into a fixed tube with impermeable gas.

 For constant coefficients $(\al=\beta=0)$  with large initial data,  Kazhikhov and Shelukhin
\cite{9} first obtained   the global existence of solutions  under boundary conditions \eqref{1.77}. From then on,
significant progress has been made on the mathematical aspect of the
initial
boundary value problems, see \cite{1,2,3,15,17,18} and the references therein.
Moreover, much effort has been made to generalize this approach to other cases.   Motivated by the fact that in the case of isentropic flow a temperature dependence on the viscosity translates
into a density dependence, there is a body of literature (see \cite{hs1,hsj1,HuangShi,3,ssyy11,28}   and the references therein) studying  the case that $\mu$  is
independent of $\te$, and heat conductivity is allowed to depend on temperature in a special
way with a positive lower bound and balanced with corresponding constitution relations.

 Kawohl \cite{ssyy11}, Jiang \cite{ssyy22,ssyy33}  and Wang \cite{ssyy44} established the
global existence of smooth solutions for \eqref{1.1}--\eqref{1.3}, \eqref{1.6} with boundary condition of either \eqref{1.7} or
\eqref{1.77} under the assumption $\mu(v) \geq \mu_0> 0$ for any $v>0$ and $\kappa$ may depend on both density and
temperature.    However, it should be mentioned here that the methods used there relies heavily on the non-degeneracy of both the viscosity $\mu$ and the heat conductivity $\ka$  and  cannot be applied directly to  the degenerate and nonlinear case  ($\al\geq 0,\b>0$).   Under the assumption that $\al=0$ and $\beta\in(0,3/2),$ Jenssen-Karper \cite{24}  proved the global existence of a weak solution to \eqref{1.1}--\eqref{1.7}.
 Later, for $\al=0$ and $\beta\in(0,\infty),$ Pan-Zhang\cite{28} obtain the global strong solutions. In \cite{24,28}, they only consider the case of non-slip and heat insulated boundary conditions.
  Recently, for the case of   stress-free and heat insulated boundary condition,  Duan-Guo-Zhu \cite{7} obtain the global strong solutions of \eqref{1.1}-\eqref{1.7}  under the condition that \be \la{h1.8}( v_0,u_0,\te_0)\in H^1\times H^2\times H^2. \ee

  In fact,  one of the main aims of this paper is  to prove  the existence and uniqueness of global strong solutions to  \eqref{1.1}-\eqref{1.7}  for $\al\geq 0 $ and $\b>0$ with  the  conditions on the initial data: \bnn ( v_0,u_0,\te_0)\in H^1 , \enn which is
   similar as those of \cite{9}. 
 Then  we state
  our main result   as follows.
 \begin{theorem}\la{thm1.1} Suppose that \be  \la{1.8}\al\geq 0,\quad \beta> 0,\ee
 and that the initial data $ ( v_0,u_0,\te_0)$   satisfies
  \be\la{h1.10} ( v_0,u_0,\te_0)\in   H^1 (0,1),\ee  and \be\la{1.9}
\inf_{x\in (0,1)}v_0(x)>0, \quad \inf_{x\in (0,1)}\theta_0(x)>0. \ee
Then, the initial-boundary-value problem \eqref{1.1}-\eqref{1.7} has a unique strong solution $(v,u,\te)$ such that for each fixed $ T>0 $,
 \be
 \begin{cases} \la{1.10}  v ,\,u,\,\theta \in L^\infty(0,T;H^1(0,1)),\\ v_t\in
  L^\infty(0,T;L^2(0,1))\cap L^2(0,T;H^1(0,1)), \\ u_t,\,\theta_t,\,u_{xx},\,\te_{xx} \in
  L^2((0,1)\times(0,T)),\end{cases}\ee
  and
  \be \inf_{(x,t)\in(0,1)\times(0,T)} v(x,t)\ge C^{-1} ,\quad \inf_{(x,t)\in(0,1)\times(0,T)}\te(x,t)\geq C^{-1},\ee
  where $C $ is a positive constant depending only on the data and $T.$
  \end{theorem}

  A few remarks are in order.
\begin{remark} Our Theorem \ref{thm1.1} can be regarded as a natural generalization of the Kazhikhov-Shelukhin's result (\cite{9})  for the constant heat conductivity case to the degenerate and nonlinear one under stress-free and thermally insulated boundary conditions.\end{remark}

  \begin{remark}Our result  improves  Duan-Guo-Zhu's result \cite{7} where they need the initial data satisfy  \eqref{h1.8}   which are   stronger than \eqref{h1.10}.  \end{remark}

We now   comment on the analysis of this paper. After modifying slightly the method due to Kazhikhov-Shelukhin \cite{9}, we obtain a key representation of $v$ (see \eqref{eq2.3}) which can be  used to obtain directly not only the lower bound of $v$ (see \eqref{eq3.1a}) but also a pointwise estimate between $v$ and $\te$ (see \eqref{rem1}). A direct consequence of this pointwise estimate between $v$ and $\te$  (see \eqref{rem1}) is the bound on $L^\infty(0,T;L^1(0,1))$-norm of $v$ (see \eqref{eqhh}) which play an important role in getting the upper bound  of $v$ but cannot be obtained directly from \eqref{1.1} due to the stress-free boundary condition \eqref{1.7}. Next,   we multiply the momentum equation \eqref{1.2} by  $ (\frac{\mu}{v}u_x-P)_x$ and make full use of the stress-free boundary condition to find that the $L^2((0,1)\times (0,T))$-norm of $u_{xx}$ can be bounded by the $L^2((0,1)\times (0,T))$-norm of $\te^\b\te_x$ (see \eqref{nna1}) which indeed can be obtained by multiplying the equation of $\te$ (see \eqref{eq3.5}) by $\te^{1+\b}$ and using Gronwall's inequality (see \eqref{ppp}).  Once we get the bounds on the $L^2((0,1)\times (0,T))$-norm of both $u_{xx}$ and $u_t$ (see \eqref{yyy1}), the desired estimates on $\te_t$ and $\te_{xx}$ can be  obtained by standard method (see \eqref{6.0}). The whole procedure will be carried out in the next section.

\section{Proof of \thmref{thm1.1}}
We first state   the following   local existence result which can be proved by using the principle of compressed mappings (c.f. \cite{10,13,tan}).

\begin{lemma} Let \eqref{1.8}-\eqref{1.9} hold. Then there exists some $T>0$ such that  the initial-boundary-value problem \eqref{1.1}-\eqref{1.7} has a unique strong solution $(v,u,\te)$ satisfying\bnn
 \begin{cases}   v , \,u, \, \theta  \in L^\infty(0,T;H^1(0,1)),\\ v_t\in
  L^\infty(0,T;L^2(0,1))\cap  L^2(0,T;H^1(0,1)), \\ u_t,\,\theta_t,\,v_{xt},\,u_{xx},\,\theta_{xx} \,\in
  L^2((0,1)\times(0,T)).\end{cases}\enn \end{lemma}

Then, the proof of   Theorem \ref{thm1.1} is based on the use of a priori   estimates (see   \eqref{eee}, \eqref{yyy}, \eqref{yyy1}, and \eqref{6.0} below) the constants in
which depend only on the data of the problem.
The estimates make it possible to continue the local solution  to the whole interval $[0,\infty)$.
Without loss of generality, we assume that $\ti\mu= \ti\ka=R=c_v=1 $. 

 Next, we derive the following representation of $v$ which is essential in obtaining the time-depending upper and lower bounds of $v$.
\begin{lemma}\la{lemma20}
We have the following expression of v
\be\ba\la{eq2.3} v(x,t)=B_0(x)D_1(x,t)D_2(x,t)\left\{1+\frac{k(\al)}{B_0(x)}\int_0^t\frac{\te(x,
\tau)}{ D_1(x,\tau) D_2(x,\tau)}d\tau\right\} ,          \ea\ee
where
\be\ba\la{eq2.5}  B_0(x)=\begin{cases}\exp\left(\ln v_0(x)-\frac{1}{\al v_0(x)^\al}\right), &\mbox{ if } \alpha>0, \\v_0(x), &\mbox{ if } \alpha= 0,\end{cases} \ea\ee
\be\ba\la{eq2.4} D_1(x,t)=\exp\left\{k(\al)\int_0^x\left(u(y,t)-u_0(y)\right)dy\right\},  \ea\ee \be \la{eq2.4'} D_2(x,t)=\begin{cases}\exp\left\{\frac{1}{\al v(x,t)^\al}\right\}, &\mbox{ if } \alpha>0,\\1, &\mbox{ if } \alpha= 0, \end{cases} \ee and  \be k(\al)=\begin{cases}1, &\mbox{ if } \alpha>0,\\1/2, &\mbox{ if } \alpha= 0 . \end{cases} \ee
\end{lemma}
\pf First, it follows from \eqref{1.2} that
 \bnn u_t=\left(\frac{\mu }{v}u_x-P\right)_x, \enn
Integrating this over $(0,x)$ and using \eqref{1.7} gives
\be\ba\la{eq2.8}\left(\int_0^xudy\right)_t&=\frac{\mu }{v}u_x-P.\ea\ee
Then, on the one hand, if $\al>0,$  since $u_x=v_t$, we have
\be\ba\la{eq2.9}
\left(\int_0^xudy\right)_t=\left(\ln v-\frac{1}{\al  v^{\al}}\right)_t-\frac{\te}{v}.
\ea\ee
Integrating $\eqref{eq2.9}$ over $(0,t)$  yields
\bnn  \ln v-\frac{1}{\al v^{\al}}-\ln v_0+\frac{1}{\al v_0^\al}-\int_0^t\frac{\te}{v}d\tau =\int_0^x(u-u_0)dy,\enn
 which implies
\be\ba\la{eq2.10a} v(x,t)=B_0(x)D_1(x,t)D_2(x,t) \exp\left\{k(\al)\int_0^t\frac{\te}{v}(x,\tau)d\tau\right\},          \ea\ee
 with $D_1(x,t),$ $D_2(x,t)$ and $B_0(x)$ as in   \eqref{eq2.5}-\eqref{eq2.4'}   respectively.
On the other hand, if $\al=0$,
it follows from \eqref{eq2.8} that
\bnn\left(\int_0^x udy\right)_t=2(\ln v)_t-\frac{\te}{v}.\enn
Integrating this over $(0,t)$ leads to
\bnn\la{eq2.12} 2\ln v-2\ln v_0 =\int_0^x(u-u_0)dy+\int_0^t\frac{\te}{v}d\tau,\enn
which shows \eqref{eq2.10a} still holds for $\al=0$.
Finally, denoting
\bnn\ba  g(x,t)=k(\al)\int_0^t\frac{\te}{v}(x,\tau)d\tau,\ea\enn
we have by  \eqref{eq2.10a}
\bnn\ba g_t=\frac{k(\al)\te(x,t)}{v(x,t)}=\frac{k(\al)\te(x,t)}{ B_0(x)D_1(x,t) D_2(x,t)\exp \{g\}} ,\ea\enn
which gives
\bnn\ba  \exp\{g\}  =1+\frac{k(\al)}{B_0(x)}\int_0^t\frac{\te(x,
\tau)}{ D_1(x,\tau) D_2(x,\tau)}d\tau .\ea\enn Putting this into \eqref{eq2.10a} yields \eqref{eq2.3} and finishes the proof of Lemma \ref{lemma20}.
\thatsall

   With Lemma  \ref{lemma20} at hand,   we are in a position to prove  the   lower bounds of $v$ and $\te$.
\begin{lemma}\la{lemma30}
It holds
\be\ba\la{eq3.1} \min_{(x,t)\in[0,1]\times[0,T]}v(x,t)\geq C^{-1},\quad \min_{(x,t)\in[0,1]\times[0,T]}\te(x,t)\geq C^{-1},   \ea\ee
where (and in what follows) $C$ denotes generic positive constant depending only on $\b,\al, T, \|(v_0,u_0,\theta_0)\|_{H^1(0,1)},
 \inf\limits_{x\in (0,1)}v_0(x), $ and $ \inf\limits_{x\in (0,1)}\theta_0(x).$
\end{lemma}
\pf
First, integrating \eqref{1.3} over $(0,1)$ and using \eqref{1.7} immediately leads to \be\ba\la{eq2.2}\int_0^1\left(\te+\frac{1}{2}u^2\right)(x,t)dx=\int_0^1\left(\te+\frac{1}{2}u^2\right)(x,0)dx , \ea\ee which in particular gives
\bnn \left|\int_0^xudy\right| \leq\int_0^1|u|dy
\leq\left(\int_0^1u^2dy\right)^{\frac{1}{2}}  \leq C.\enn  Combining this  with \eqref{eq2.4} implies
\be\ba\la{eq3.2} C^{-1}\le  D_1(x,t)\le  C, \ea\ee which together with \eqref{eq2.3} yields that for any $(x,t)\in[0,1]\times[0,T],$
\be\ba\la{eq3.1a} v(x,t)\geq B_0(x)D_1(x,t)D_2(x,t)\ge C^{-1},          \ea\ee due to   \bnn\la{eq3.1b}  D_2(x,t)\ge 1,\quad C^{-1}\leq B_0(x)\leq C.\enn
  Finally, we  rewrite \eqref{1.3} as
\be\ba\la{eq3.5}
\te_t+\frac{\te}{v}u_x=\left(\frac{\te^\b \te_x}{v}\right)_x+\frac{\mu u_x^2}{v}.\ea\ee
For $r>2,$  multiplying the above equality by $\te^{-r}$ and integrating the resultant equality over $(0,1)$ yields that
\be\ba\la{eq3.6} & \frac{1}{r-1}\frac{d}{dt}\int_0^1\left( {\te}^{-1}\right)^{r-1}dx+\int_0^1\frac{\mu u_x^2}{v\te^r}dx+r\int_0^1\frac{\te^\b\te_x^2}{v\te^{r+1}}dx\\&=\int_0^1\frac{u_x}{v\te^{r-1}}dx\\&\leq
\frac{1}{2}\int_0^1\frac{\mu u_x^2}{v\te^r}dx+\frac{1}{2}\int_0^1\frac{1}{\mu v\te^{r-2}}dx\\&\leq\frac{1}{2}\int_0^1\frac{\mu u_x^2}{v\te^r}dx+C\left\|  \te^{-1} \right\|^{r-2}_{L^{r-1}} ,
 \ea\ee
where in the second inequality we have used  $\mu v=v+v^{1-\al}>v\geq C^{-1}$. Combining \eqref{eq3.6} with  Gronwall's inequality yields
\bnn\ba \sup_{0\le t\le T}\left\|  \te^{-1}(\cdot,t) \right\|_{L^{r-1}}\leq C,\ea\enn with $C$
 independent of $r.$ Letting $r\rightarrow\infty$  proves   the second inequality of \eqref{eq3.1}. Thus, the proof of
\lemref{lemma30} is finished. \thatsall

\begin{lemma} \la{lemma40}  There exists a positive constant $C$ such that  for each $(x,t)\in [0,1]\times[0,T],$
\be\ba\la{eee} C^{-1}\leq v(x,t)\leq   C.\ea\ee
\end{lemma}
\pf  First, it follows from \eqref{eq3.1} and \eqref{eq2.4'} that for any $(x,t)\in[0,1]\times[0,T],$\bnn 1\leq D_2(x,t)\leq C, \enn
which together with \eqref{eq2.3} and \eqref{eq3.2}   yields that for any $(x,t)\in[0,1]\times[0,T],$
\be\ba\la{rem1}C^{-1} \leq v(x,t)\leq C+C\int_0^t\te(x,\tau) d\tau.\ea\ee
Integrating this with respect to $x$ over $(0,1) $ and using \eqref{eq2.2} leads to \be\ba\la{eqhh} \sup_{0\le t\le T}\int_0^1 v(x,t)dx\leq C.\ea\ee
Next, for $\eta \in (0,1) $ and $\ve \in (0,1),$
 integrating \eqref{eq3.5} multiplied by $\te^{-\eta}$ over $ (0,1)\times(0,T),$ we get by   \eqref{eq2.2} and \eqref{eq3.1}
\be\ba\la{eqmm}
&\int_0^T\int_0^1\frac{\eta\te^\b\te_x^2}{v\te^{\eta+1}}dxdt+\int_0^T\int_0^1\frac{\mu u_x^2}{v\te^{\eta}}dxdt\\&
=\frac{1}{1-\eta} \int_0^1\te^{1-\eta}dx-\frac{1}{1-\eta}\int_0^1\te_0^{1-\eta}dx +\int_0^T\int_0^1\frac{u_x}{\te^{\eta-1 }v}dxdt \\&\leq C(\eta)+\frac{1}{2}\int_0^T\int_0^1\frac{\mu u_x^2}{v\te^\eta}dxdt+C\int_0^T\int_0^1\frac{1}{\mu v\te^{\eta-2}}dxdt \\&\leq
C(\eta)+\frac{1}{2}\int_0^T\int_0^1\frac{\mu u_x^2}{v\te^\eta}dxdt+C\int_0^T\max_{x\in[0,1]}\te^{1-\eta}dt\\&\leq
C(\eta,\ve)+\frac{1}{2}\int_0^T\int_0^1\frac{\mu u_x^2}{v\te^\eta}dxdt+\ve\int_0^T\max_{x\in[0,1]}\te dt,\ea\ee
where in the first inequality we have used
\bnn \int_0^1\te^{1-\eta}dx\leq C .  \enn
  Finally, using \eqref{eqhh},  we obtain that   for $ \eta=\min\{1,\b\}/2,$
\bnn\ba\int_0^T\max_{x\in[0,1]}\te dt&\leq C+C\int_0^T\int_0^1|\te_x|dxdt \\&\leq C+C\int_0^T\int_0^1\frac{  \te^\b\te_x^2}{v\te^{1+\eta}}dxdt+C\int_0^T\int_0^1\frac{v\te^{1+\eta}}{ \te^\b}dxdt\\&\leq C+C\int_0^T\int_0^1\frac{  \te^\b\te_x^2}{v\te^{1+\eta}}dxdt+C\int_0^T\max_{x\in[0,1]}\te^{1+\eta-\b}dt
 \\&\leq C +C\int_0^T\int_0^1\frac{  \te^\b\te_x^2}{v\te^{1+\eta}}dxdt+\frac{1}{2}\int_0^T\max_{x\in[0,1]}\te dt,\ea\enn
   which together with \eqref{eqmm} yields that \be \la{tem}\int_0^T\max_{x\in[0,1]}\te dt\le C,\ee and that  for $\eta \in (0,1) $ \be\ba\la{eqaa}\int_0^T\int_0^1  \te^{\b-1-\eta}\te_x^2 dxdt  \leq C(\eta). \ea\ee
 Combining   \eqref{rem1} with  \eqref{tem}
  finishes the proof of \lemref{lemma40}. \thatsall

\begin{lemma}\la{lemma50}
There exists a positive constant $C$ such that
\be\la{yyy}\sup_{0\leq t\leq T}\int_0^1 v_x^2 dx \leq C.  \ee
\end{lemma}
\pf First, we rewrite \eqref{1.2} as\bnn\ba\left(\frac{\mu v_x}{v}-u\right)_t=\left(\frac{\te}{v}\right)_x.\ea\enn
Multiplying the above equality by $\frac{\mu v_x}{v}-u$ and integrating  it over $(0,1)\times(0,T)$ gives
\be\ba\la{aaa}&
\frac{1}{2} \int_0^1\left(\frac{\mu v_x}{v}-u\right)^2dx-\frac{1}{2}\int_0^1\left(\frac{\mu v_x}{v}-u\right)^2(x,0)dx+\int_0^T\int_0^1\frac{\mu \te v_x^2}{v^3}dxdt\\&=
 \int_0^T\int_0^1\frac{\te v_x u}{v^2}dxdt +\int_0^T\int_0^1\frac{\te_x}{v}\left(\frac{\mu v_x}{v}-u\right)dxdt
 \triangleq I_1+I_2.
\ea\ee
Then, on the one hand, Cauchy's inequality, \eqref{eq2.2}, \eqref{eee}, and \eqref{tem}  lead to
\be\ba\la{bbb} |I_1|&\leq\frac{1}{2}\int_0^T\int_0^1\frac{\te u^2}{v\mu }dxdt+\frac{1}{2}\int_0^T\int_0^1\frac{\mu \te v_x^2}{v^3}dxdt\\&\leq C\int_0^T\max_{x\in[0,1]}\te dt+\frac{1}{2}\int_0^T\int_0^1\frac{\mu \te v_x^2}{v^3}dxdt\\&\leq C+\frac{1}{2}\int_0^T\int_0^1\frac{\mu \te v_x^2}{v^3}dxdt.
\ea\ee
On the other hand, we deduce  from  \eqref{eee}, \eqref{eqaa}, and \eqref{eq3.1} that for $ \eta=\min\{1,\b\}/2,$
\be\ba\la{ccc} |I_2|&\leq C\int_0^T\int_0^1 \te^{\b-1-\eta} \te_x^2 dxdt+C\int_0^T\int_0^1  \te^{1+\eta-\b}  \left(\frac{\mu v_x}{v}-u\right)^2dxdt\\&
\leq C+C\int_0^T\max_{x\in[0,1]}\te^2\int_0^1\left(\frac{\mu v_x}{v}-u\right)^2dxdt.
 \ea\ee
Next, it follows from \eqref{eq2.2}  and \eqref{eee}  that  for $ \eta=\min\{1,\b\}/2,$
 \be\ba\la{sy3}&\int_0^T\max_{x\in[0,1]}(\te^{1+\b }+\te^2)dt\\& \le C \int_0^T\max_{x\in[0,1]}\te^{2+\b-\eta}dt\\&\leq C
 +C\int_0^T\left(\max_{x\in[0,1]}\left|\te^{\frac{2+\b-\eta}{2}}(x,t) -(\int_0^1\te dx)^{\frac{2+\b-\eta}{2}}\right|\right)^2dt\\&
 \leq C+C\int_0^T\left(\int_0^1\te^{\frac{\b-\eta}{2}}\left|\te_x\right|dx\right)^2dt \\&\leq C+ C\int_0^T\left(\int_0^1\frac{ \te^\b\te_x^2}{v\te^{\eta+1}}dx\right)\left(\int_0^1  {v\te}  dx\right)dt\\&
 \leq C+C\int_0^T\int_0^1\frac{ \te^\b\te_x^2}{v\te^{\eta+1}}dxdt
\\&\leq C, \ea\ee where in the last inequality we have used \eqref{eqaa}.
Finally, adding \eqref{bbb} and \eqref{ccc} to \eqref{aaa}, we obtain after using Gronwall's inequality and \eqref{sy3} that
\bnn\ba \sup_{0\le t\le T}\int_0^1\left(\frac{\mu v_x}{v}-u\right)^2dx+\int_0^T\int_0^1\frac{\mu \te v_x^2}{v^3}dxdt\leq C,  \ea\enn
which together with \eqref{eq2.2}  and \eqref{eee} gives
\be\ba\la{ggg}\sup_{0\le t\le T}\int_0^1v_x^2dx&\leq C+C\sup_{0\le t\le T}\int_0^1\left(\frac{\mu v_x}{v}-u\right)^2dx\\&\leq C.\ea\ee
 The proof of \lemref{lemma50} is finished. \thatsall

\begin{lemma}\la{lemma5a}
There exists a positive constant $C$ such that
\be \la{yyy1}\sup_{0\leq t\leq T}\int_0^1 u_x^2 dx+\int_0^T\int_0^1(u_t^2+u_{xx}^2)dxdt\leq C.  \ee
\end{lemma}
\pf First, integrating \eqref{1.2} multiplied by $ (\frac{\mu}{v}u_x-P)_x$ over $(0,1) ,$  we  obtain by   integration by parts and \eqref{eq3.5}  that
\bnn\ba & \int_0^1  \left(\frac{\mu}{v}u_x-P\right)_x^2dx\\&=-\int_0^1  \left(\frac{\mu}{v}u_x-P\right)u_{tx}dx\\&=-\frac12\int_0^1 \frac{\mu}{v}\left(u^2_x\right)_tdx +\left(\int_0^1  P u_{ x}dx\right)_t-\int_0^1  P_t u_{ x}dx\\&=-\frac12\left(\int_0^1 \frac{\mu}{v}u^2_xdx\right)_t+\frac12 \int_0^1 \left(\frac{\mu}{v}\right)_v' u^3_xdx +\left(\int_0^1  P u_{ x}dx\right)_t\\&\quad+2\int_0^1  \frac{\te}{v^2} u^2_{ x}dx +\int_0^1 \frac{\te^\b\te_x}{v}\left(\frac{u_x}{v}\right)_xdx -\int_0^1 \frac{\mu u_x^3}{v^2}  dx,\ea\enn
which in particular gives
\be\la{nna1}\ba &  \left(\int_0^1\left(\frac{\mu}{2v}u^2_x-Pu_x\right)dx\right)_t+\int_0^1 \left(\frac{\mu^2}{v^2}u_{xx}^2+ \frac{\te_x^2}{v^2}\right)dx\\& \le C\int_0^1|u_{xx}|\left(|v_x||u_x|+|v_x|\te+ \te^\b |\te_x|\right)dx +C\int_0^1|\te_x| |v_x|\te dx\\&\quad+C\int_0^1|\te_x| \te^\b |u_x||v_x|dx+C\int_0^1 v_x^2\left(u_x^2+\te^2\right)dx+C\int_0^1\left(|u_x|^3+u_x^2\te\right)dx\\& \le \frac14\int_0^1  \frac{\mu^2}{v^2}u_{xx}^2dx+ C\int_0^1 \te^{2\b}  \te_x^2dx +C\left(\int_0^1 u_x^2 dx\right)^2\\&\quad +C \max_{x\in[0,1]} \left(u_x^2+ \te^2 \right)\left(1+\int_0^1 v_x^2dx+\int_0^1 \te dx\right)\\& \le \frac12\int_0^1  \frac{\mu^2}{v^2}u_{xx}^2dx+ C_1\int_0^1 \frac{\te^{2\b}\te_x^2}{v}dx +C\max_{x\in[0,1]}\te^2   +C\left(\int_0^1 u_x^2 dx\right)^2+C,\ea\ee
where in the last inequality we have used \eqref{ggg} and \be\ba\la{heq2} \max_{x\in[0,1]} u_x^2 \le C(\ve)\int_0^1 u_x^2dx+\ve \int_0^1 u_{xx}^2dx,\ea\ee for $\ve>0$ small enough.
Then,  integrating \eqref{eq3.5} over $(0,1)\times(0,T)$ yields that \bnn\ba \int_0^T\int_0^1 \frac{\mu u_x^2}{v}dxdt &=\int_0^1\te dx-\int_0^1\te_0dx +\int_0^T\int_0^1 \frac{\te}{v}u_xdxdt\\&\le C+ \frac12\int_0^T\int_0^1 \frac{\mu u_x^2}{v}dxdt+C\int_0^T\int_0^1 \frac{\te^2}{\mu v} dxdt\\&\le C+\frac12\int_0^T\int_0^1 \frac{\mu u_x^2}{v}dxdt+C\int_0^T\max_{x\in [0,1]}\te(x,t)dt\ea\enn
which together with \eqref{tem} and \eqref{eee} gives
\be\ba\la{sy6} \int_0^T\int_0^1u_x^2dxdt\leq C.\ea\ee
Next, integrating \eqref{eq3.5} multiplied by $\te^{1+\b}$ over $ (0,1) $ leads to
\be\ba\la{nnn}
&\frac{1}{2+\b}(\int_0^1\te^{2+\b}dx)_t+ (1+\b)\int_0^1\frac{\te^{2\b}\te_x^2}{v}dx \\&=
 - \int_0^1\frac{\te^{2+\b}u_x}{v}dx + \int_0^1\frac{\mu \te^{1+\b} u_x^2}{v}dx \\&\leq
 C \int_0^1\te^{3+\b}dx+C  \int_0^1 \te^{1+\b}  u_x^2dx   \\&\leq
C \max_{x\in[0,1]}\te\int_0^1\te^{2+\b}dx +C \max_{x\in[0,1]}\te^{1+\b}\int_0^1   u_x^2dx  .
\ea\ee

Choosing $C_2\ge C_1+1$ suitably large such that
\bnn C_2\te^{2+\b}+\mu v^{-1}u_x^2\ge 4(2+\b)\te v^{-1}|u_x|,\enn
 adding \eqref{nnn} multiplied by $C_2$ to  \eqref{nna1}, and choosing $\varepsilon$  sufficiently small, we obtain from Gronwall's inequality, \eqref{sy6},   and \eqref{sy3} that
\be\ba\la{ppp}\sup_{0\leq t\leq T}\int_0^1\left(\te^{2+\b}+u_x^2\right)dx+\int_0^T\int_0^1 u_{xx}^2 dxdt
+\int_0^T\int_0^1 {\te^{2\b}\te_x^2} dxdt\leq C.  \ea\ee
Finally, rewriting \eqref{1.2} as
\bnn\ba\la{eq20} u_t=\frac{\mu u_{xx}}{v}+\left(\frac{\mu }{v}\right)_v' v _xu_x-\frac{\te_x}{v}+\frac{\te v_x}{v^2},\ea\enn
  we deduce from  \eqref{eee}, \eqref{ppp}, \eqref{ggg}, \eqref{sy6}, \eqref{heq2}, and \eqref{sy3} that
\bnn\ba
\int_0^T\int_0^1u_t^2dxdt&\leq C\int_0^T\int_0^1\left(u_{xx}^2+u_x^2v_x^2+\te_x^2+\te^2v_x^2\right)dxdt\\&\leq C+C\int_0^T\max_{x\in[0,1]}\te^2dt\\&\leq C,
\ea\enn
which together with \eqref{ppp} finishes the proof of Lemma \ref{lemma5a}. \thatsall

\begin{lemma}\la{lemma60}There exists a positive constant $C$ such that
\be\la{6.0} \sup_{0\le t\le T}\xix  \te_x^2dx+\int_0^T\xix \left(  \te_t^2+\te_{xx}^2\right)dxdt\le C.\ee
\end{lemma}
\pf First, multiplying \eqref{eq3.5} by $\te^{\b}\te_t$ and integrating the resultant equality over (0,1) yields
\bnn\ba&\int_0^1\te^\b\te_t^2dx+\int_0^1\frac{\te^{\b+1}u_x\te_t}{v}dx\\&=-\int_0^1\left(\frac{\te^\b\te_x}{v}\right)\left(\te^\b\te_x\right)_tdx+\int_0^1\frac{\mu u_x^2}{v}\te^\b\te_tdx\\&
=-\frac{1}{2}\frac{d}{dt}\int_0^1\frac{(\te^\b\te_x)^2}{v}dx+\frac{1}{2}\int_0^1\left(\te^\b\te_x\right)^2\left(\frac{1}{v}\right)_tdx+\int_0^1\frac{\mu u_x^2\te^\b\te_t}{v}dx\\&
=-\frac{1}{2}\frac{d}{dt}\int_0^1\frac{(\te^\b\te_x)^2}{v}dx-\frac{1}{2}\int_0^1\frac{(\te^\b\te_x)^2u_x}{v^2}dx +\int_0^1\frac{\mu u_x^2\te^\b\te_t}{v}dx,\ea\enn
which combined with the H\"{o}lder inequality,  \eqref{yyy}, and \eqref{ppp} leads to
\be\ba\la{equa00}&\frac{1}{2}\frac{d}{dt}\int_0^1\frac{\left(\te^\b\te_x\right)^2}{v}dx+\int_0^1\te^\b\te_t^2dx\\&=
-\frac{1}{2}\int_0^1\frac{\left(\te^\b\te_x\right)^2u_x}{v^2}dx+\int_0^1\frac{\mu u_x^2\te^\b\te_t}{v}dx-\int_0^1\frac{\te^{\b+1}u_x\te_t}{v}dx\\&   \leq C \int_0^1\te^{2\b} \te_x^2|u_x|dx + \frac{1}{2}\int_0^1\te^\b\te_t^2dx +C\int_0^1u_x^4\te^\b dx+C\int_0^1\te^{\b+2}u_x^2dx \\&\leq C\max_{x\in[0,1]}|u_x| \int_0^1\te^{2\b} \te_x^2dx +\frac{1}{2}\int_0^1\te^\b\te_t^2dx+
C\max_{x\in[0,1]}\left(u_x^2\te^\b+\te^{\b+2}\right)\\&\leq
\frac{1}{2}\int_0^1\te^\b\te_t^2dx+C(\int_0^1\te^{2\b}\te_x^2dx)^2 +
C\max_{x\in[0,1]}\left(u_x^4+\te^{2\b+2}\right)+C.
\ea\ee
It follows from \eqref{yyy}, \eqref{heq2}, \eqref{yyy1}, and the H\"{o}lder  inequality that
 \be\ba\la{equa11}\int_0^T\max_{x\in[0,1]}u_x^4dt&\leq C\int_0^T\int_0^1u_x^4dxdt+C\int_0^T\int_0^1|u_x^3u_{xx}|dxdt\\&\leq C\int_0^T\max_{x\in[0,1]}u_x^2\int_0^1u_x^2dxdt\\&\quad+C\int_0^T\max_{x\in[0,1]}u_x^2
 \left(\int_0^1u_x^2dx\right)^{\frac{1}{2}}\left(\int_0^1u_{xx}^2dx\right)^\frac{1}{2}dt\\&\leq C\int_0^T\int_0^1\left(u_x^2 + u_{xx}^2\right)dxdt +\frac12 \int_0^T \max_{x\in[0,1]}u_x^4dt
 \\&\leq C +\frac12 \int_0^T \max_{x\in[0,1]}u_x^4dt \ea\ee
which combined with \eqref{equa00},  \eqref{ppp}, and the Gronwall  inequality yields
\be\ba\la{equa22} \sup_{0 \le t\le T}\xix  \left(\te^\b\theta_{x}\right)^2 dx+\int_0^T\xix  \te^\b\te_t^2dxdt\le C, \ea\ee
where we have used  
 \be\ba\la{equa44}\max_{x\in[0,1]}\te^{2\b+2} &\leq C+C\left(\max_{x\in[0,1]}|\te^{\b+1}(x,t)-(\int_0^1\te dx)^{\b+1} |\right)^2
 \\&\leq
 C+C\left(\int_0^1|\te^\b\te_x|dx\right)^2\\&\leq
 C+C\int_0^1(\te^\b\te_x)^2dx.    \ea  \ee
Combining \eqref{equa22} with \eqref{equa44}  implies for all $(x,t)\in (0,1)\times(0,T)$\be\la{equa55} \te (x,t) \le C .\ee
Meanwhile, both \eqref{equa22} and \eqref{eq3.1} lead  to
\be\ba\la{equa66}\sup_{0 \le t\le T}\xix  \theta_{x}^2 dx+\int_0^T\xix  \te_t^2dxdt\le C. \ea\ee
Finally, it follows from \eqref{eq3.5} that
\bnn\ba\frac{\te^\b\te_{xx}}{v}=-\frac{\b\te^{\b-1}\te_x^2}{v}+\frac{\te^\b\te_xv_x}{v^2}-\frac{\mu u_x^2}{v}+\frac{\te u_x}{v}+\te_t,    \ea\enn
 which together with \eqref{equa66}, \eqref{eq3.1}, \eqref{equa55}, \eqref{equa11}, and \eqref{yyy} gives
\bnn\ba\int_0^T\int_0^1\te_{xx}^2dxdt&\leq C\int_0^T\int_0^1\left(\te_x^4+v_x^2\te_x^2+u_x^4+\te^2u_x^2+\te_t^2\right)dxdt\\&\leq
C+C\int_0^T\max_{x\in[0,1]}\te_x^2dt\\&\leq
C+C\int_0^T\int_0^1\te_x^2dxdt+\frac{1}{2}\int_0^T\int_0^1\te_{xx}^2dxdt.
\ea\enn
 Combining this with \eqref{equa66} gives \eqref{6.0} and   finishes  the proof of \lemref{lemma60}.\thatsall

 \end{document}